\documentclass[12pt]{iopart}
\usepackage{iopams}
\usepackage{graphicx}
\input amssym.def \input amssym

\begin{document}

\title[]{Riemann hypothesis from the Dedekind psi function}

\author{Michel Planat 
}

\address{ Institut FEMTO-ST, CNRS, 32 Avenue de
l'Observatoire,\\ F-25044 Besan\c con, France. }
\vspace*{.1cm}

\begin{abstract}
Let $\mathcal{P}$ be the set of all primes and $\psi(n)=n\prod_{n\in \mathcal{P},p|n}(1+1/p)$ be the Dedekind psi function. We show that the Riemann hypothesis is satisfied if and only if $f(n)=\psi(n)/n-e^{\gamma} \log \log n <0$ for all integers $n>n_0=30$ (D), where $\gamma \approx 0.577$ is Euler's constant.  This inequality is equivalent to Robin's inequality that is recovered from (D) by replacing $\psi(n)$ with the sum of divisor function $\sigma(n)\ge \psi(n)$ and the lower bound by $n_0=5040$. For a square free number, both arithmetical functions $\sigma$ and $\psi$ are the same. We also prove that any exception to (D) may only occur at a positive integer $n$ satisfying $\psi(m)/m<\psi(n)/n$, for any $m<n$, hence at a primorial number $N_n$ or at one its multiples smaller than  $N_{n+1}$ (Sloane sequence $A060735$). According to a Mertens theorem, all these candidate numbers are found to satisfy (D): this implies that the Riemann hypothesis is true.
\end{abstract} 


\noindent

\section{Introduction}

Riemann zeta function $\zeta(s)=\sum_{n>1}n^{-s}=\prod_{p \in \mathcal{P}}\frac{1}{1-p^{-s}}$ (where the product is taken over the set $\mathcal{P}$ of all primes) converges for $\mathcal{R}(s)>1$. It has a analytic continuation to the complex plane with a simple pole at $s=1$. The Riemann hypothesis (RH) states that non-real zeros all lie on the critical line $\mathcal{R}(s)=\frac{1}{2}$. RH has equivalent formulations, many of them having to do with the distribution of prime numbers \cite{Edwards74,Hardy79}.

Let $d(n)$ be the divisor function. There exists a remarkable parallel between the error term $\Delta(x)=\sum_{n \le x} d(n) -x (\log x + 2\gamma -1)$ in the summatory function of $d(n)$ (Dirichlet's divisor problem) and the corresponding mean-square estimates $|\zeta(\frac{1}{2})+it|$ of $\zeta(s)$ on the critical line, see \cite{Ivic85} for a review. This may explain Ramanujan's interest for highly composite numbers \cite{Rama88}. A highly composite number is a positive integer $n$ such that for any integer $m<n$, $d(m)<d(n)$, i.e. it has more divisors than any positive integer smaller than itself.

This landmark work eventually led to Robin's formulation of RH in terms of the sum of divisor function $\sigma(n)$ \cite{Robin84,Nicolas83,Sole07}. More precisely,
\begin{equation}
\mbox{RH}~\mbox{is}~\mbox{true}~\mbox{iff}~ g(n)=\frac{ \sigma(n)}{n} - e^{\gamma} \log \log n <0~ \mbox{for}~\mbox{any}~ n>5040.
\label{Robin}
\end{equation}
The numbers that do not satisfy (\ref{Robin}) are in the set $\mathcal{A}=\{2,3,4,5,6,8,9,10,12,16,18,20,24,$
$30,36,48,60,72,84,120,180,240,360,720,840,2520,5040\}$.

If RH is false, the smallest value of $n>5040$ that violates the inequality should be a superabundant number \cite{Alaoglu44,Akbary09,Erdos75}, i.e. a positive number $n$ satisfying
\begin{equation}
\frac{\sigma(m)}{m}<\frac{\sigma(n)}{n}~ \mbox{for}~\mbox{any}~ m<n.  
\end{equation}
No counterexample has been found so far. 

In the present paper, Robin's criterion is refined by replacing the sum of divisor function  by the Dedekind psi function $\psi(n)=\prod_{p\in \mathcal{P}, p|n} (1+\frac{1}{p})$ \footnote{The Dedekind function $\psi(n)$ should not be confused with the second Chebyshev function $\psi_T(n)=\sum_{p^k\le n}\log p.$}. Since $\psi(n) \le \sigma(n)$, with equality when $n$ is free of square,  we establish the refined criterion
\begin{equation}
\mbox{RH}~\mbox{is}~\mbox{true}~\mbox{iff}~ f(n)= \frac{\psi(n)}{n} - e^{\gamma} \log \log n <0 ~ \mbox{for}~\mbox{any}~ n>30.
\label{Dedekind}
\end{equation}
The numbers that do not satisfy (\ref{Dedekind}) are in the set $\mathcal{B}=\{2,3,4,5,6,8,10,12,18,30\}$.

A number that would possibly violate (\ref{Dedekind}) should satisfy
\begin{equation}
\frac{\psi(m)}{m}<\frac{\psi(n)}{n}~ \mbox{for}~\mbox{any}~  m<n, 
\label{newSuperab} 
\end{equation}
and be a primorial number $N_n=\prod_{i=1}^n p_i$ (the product of the first $n$ primes) or one its multiples smaller than  $N_{n+1}$ (Sloane sequence $A060735$).

According to a Mertens theorem \cite{Hardy79}, one has $\lim_{n\rightarrow \infty} \frac{\psi(N_n)}{N_n}=\frac{e^{\gamma}}{\zeta(2)} log(p_n)$ and we show that none of the numbers larger than $30$ in the sequence $A060735$ can  violate (\ref{Dedekind}). As a result, RH may only be true.

In the next section, we provide the proof of (\ref{Dedekind}) and compare it with the Robin's criterion (\ref{Robin}). Furthermore, we investigate the structure of numbers satisfying (\ref{newSuperab}) and justify why they fail to provide counterexamples to RH.

Originally, Dedekind introduced $\psi(n)$ as the index of the congruence subgroup $\Gamma_0(n)$ in the modular group (see \cite{Schoeneberg74}, p. 79). In our recent work, the Dedekind psi function plays a role for understanding the commutation relations of quantum observables within the discrete Heisenberg/Pauli group \cite{Planat10}. In particular, it counts the cardinality of the projective line $\mathbb{P}_1 (\mathbb{Z}_n)$ of the lattice $\mathbb{Z}_n \times \mathbb{Z}_n$. The relevance of the Dedekind psi function $\psi(n)$ in the context of RH is novel. For other works aiming at a refinement of Robin's inequality, we mention \cite{Lagarias02}, \cite{Sole07} and \cite{Kupershmidt09}.

\section{A proof of the refined Robin's inequality}

Let us compute the sequence $\mathcal{S}$ of all positive numbers satisfying (\ref{newSuperab}). For $2<n<10^5$, one obtains $S=\{N_1=2,4,N_2=6,12,18,24,N_3=30,60,90,120,150,180,N_4=210,420,630,840,1050,1260,1470,1680,1890,2100,N_5=2310,4620,6930,9240\}$,
which are the first terms of Sloane sequence $A060735$, consisting of the primorials $N_n$ and their multiples up to the next primorial $N_{n+1}$. 

It is straightforward to check that about half of the numbers in $S$ are not superabundant (compare to Sloane sequence $A004394$).

Based on calculations performed on the numbers in the finite sequence $S$, we are led to three properties that the infinite sequence $A060735$ should satisfy





{\bf Proposition 1}: For any $l>1$ such that $N_n<lN_n<N_{n+1}$ one has $f(lN_n)< f(N_n)$.

{\it Proof}: 
One may use $\tilde{f}(n)=\frac{\psi(n)}{n \log \log n}$ instead of $f(n)$ to simplify the proof.

When $l$ is prime, one observes that $lN_n=p_1 p_2 \cdots l^2\cdots p_n$ for some $p_j=l$. The corresponding Dedekind psi function is evaluated as 

$$\psi(l N_n)=(p_j^2+p_j)\prod_{i\ne j}\psi(p_i)=l \psi(N_n).$$

Then, with $$\tilde{f}(lN_n)=\frac{\psi(lN_n)}{lN_n\log\log(lN_n)}=\frac{\psi(N_n)}{N_n\log\log(lN_n)},$$ one concludes that
$$\frac{\tilde{f}(lN_n)}{\tilde{f}(N_n)}=\frac{\log\log N_n}{\log\log(lN_n)}<1 ~\mbox{in}~\mbox{the}~\mbox{required}~\mbox{range}~1<l<p_{n+1}.$$
 
When $l$ is not prime, a similar calculation is performed by decomposing $l$ into a product of primes and by using the multiplicative property of $\psi(n)$.

{\bf Proposition 2}: Given $l \ge 1$, for any $m$ such that $l N_n<m < (l+1)N_n<N_{n+1}$ one has  $f(m)<f(N_n)$. 

{\it Proof}: This proposition is proved by using inequality (\ref {newSuperab}) at $n=(l+1)N_n$

$$\frac{\psi(m)}{m}<\frac{\psi[(l+1)N_n]}{(l+1)N_n}~\mbox{for}~\mbox{any}~ m<(l+1)N_n$$

and the relation $\psi[(l+1)N_n]=(l+1)\psi(N_n)$. As a result

$$\tilde{f}(m)<\frac{\psi(N_n)}{N_n \log \log m}=\frac{\log \log N_n}{\log \log m}\tilde{f}(N_n)<\tilde{f}(N_n)~\mbox{since}~m>N_n.$$


{\bf Proposition 3}: There exist {\it infinitely many} prime numbers $p_n$ such that $f(N_{n+1})>f(N_n)$.

{\it Remarks on the proposition 3:} Based on experimental evidence in table 1, one would expect that $f(N_{n+1})<f(N_n)$ and that the proposition 3 is false.

Similarly, one would expect that $\theta(p_n)<p_n$ for any $n$. For instance it is known \cite{Schoenfeld76} that 

$$ \theta(n)<n ~\mbox{for}~0<n\le 10^{11}.$$

\begin{table}[ht]
\begin{center}
\small
\begin{tabular}{|r|r|r|r|r|}
\hline
$n$ & $10$ & $10^3$ & $10^5$ & $10^7$\\
\hline
$\frac{\theta(p_n)}{p_n}$ & $0.779$ & $0.986$& $0.99905$ & $0.999958$\\
$\frac{\tilde{f}(N_{n+1})}{\tilde{f}(N_n)}$& $0.987$ & $0.9999980$& $0.99999999921$ & $0.99999999999975$\\
\hline
$\frac{k_n \log k_n}{p_{n+1}\log p_{n+1}}$ & $0.938$ & $1.00378$ & $1.000447$ & $1.0000423$\\
\hline
\end{tabular}
\label{table1}
\caption{An excerpt of values of $\theta(p_n)/p_n$ and $\frac{\tilde{f}(N_{n+1})}{\tilde{f}(N_n)}$ versus the number of primes in the primorial $N_n$.}
\end{center}
\end{table}

Many oscillating functions were studied in the context of the prime number theorem. It was believed in the past that, for any real number $x$, the function $\Delta_1(x)=\pi(x)-\mbox{li}(x)$ (where \mbox{li}(x) is the logarithmic integral) is always negative. However, J~E. Littlewood has shown that $\Delta_1(x)$ changes sign infinitely often at some large values $x>x_0$ \cite{Littlewood1914}. The smallest value $x_0$ such that for the first time $\pi(x_0)>\mbox{li}(x_0)$ holds is called the Skewes number. The lowest present day value of the Skewes number is around $10^{316}$.

In what concerns the function $\Delta_4(x)=\theta(x)-x$, according to theorem 1 in \cite{Kaczorowski85} (see also \cite{Nicolas10}, Lemme 10.1), there exists a positive constant $c_4$ such that for sufficiently large $T$, the number of sign changes of $\Delta_4(x)$ in the interval $[2,T]$ is 
$$V_4(T) \ge c_4 \log T.$$ 

{\it Justification of proposition 3}

According to theorem 34 in \cite{Ingham90} (also used in \cite{Nicolas10})

$$\theta(x)-x=\Omega_{\pm}(x^{1/2}\log_3 x)~\mbox{when}~x \rightarrow \infty,$$

where $\log_3 x = \log \log \log x$. The omega notation means that there exist infinitely many real numbers $x$ satisfying
\begin{eqnarray}
&\theta(x)\ge x + \frac{1}{2}\sqrt x \log_3 x =k_x,\\
& \mbox{and}~\theta(x)\le x - \frac{1}{2}\sqrt x \log_3 x \nonumber 
\label{little}
\end{eqnarray}
If $x=p_n$, for some $n$ then
\begin{equation}
\theta(p_n)\ge p_n +\frac{1}{2}\sqrt p_n \log_3 p_n=k_n.
\label{eq}
\end{equation}
Otherwise, let us denote $p_n$ the first prime number preceeding $x$, one has
$$\theta (p_n)= \theta(x) \ge k_x\ge k_n,$$
that is similar to (\ref{eq}).

At a primorial $n=N_n$, $\psi(N_n)=\prod _{i=1}^n (1+p_i)$ so that $\psi(N_{n+1})=(1+p_{n+1})\psi(p_n)$. One would like to show that there are infinitely many prime numbers $p_n>2$ such that 
%
$$\frac{\tilde{f}(N_{n+1})}{\tilde{f}(N_n)}=(1+\frac{1}{p_{n+1}})\frac{\log \theta(p_n)}{\log \theta(p_{n+1})}=\frac{1+\frac{1}{p_{n+1}}}{1+\log(1+\frac{ \log p_{n+1}}{\theta(p_n)})/\log \theta(p_n)}>1.$$
%
By contradiction, let us assume that the reverse inequality holds for those prime numbers $x=p_n$ satisfying (\ref{little})
\begin{eqnarray}
&\frac{\log \theta(p_n)}{p_{n+1}}<\log(1+\frac{\log p_{n+1}}{\theta(p_n)})  \nonumber \\ 
& \mbox{with}~ \theta(p_n)\ge p_n + \frac{1}{2}\sqrt p_n \log_3 p_n. \nonumber \\
\label{test}
\end{eqnarray}
Taking the development of the logarithm in the first equation of  (\ref{test}) one obtains $\frac{\log k_{p_n}}{p_{n+1}}<\frac{\log p_{n+1}}{k_{p_n}}$, that is 

\begin{equation}
k_{p_n} \log k_{p_n}< p_{n+1} \log p_{n+1} ~\mbox{for}~k_{p_n}=p_n +\frac{1}{2}\sqrt p_n \log_3 p_n.
\label{ineq}
\end{equation} 

The inequality (\ref{ineq}) contradicts the calculations performed in table 1 for $10<n<10^7$. We conclude that our proposition 3 is correct in a finite range of $p_n$'s. In addition, since there are infinitely many prime numbers satisfying (\ref{little}), proposition 3 is also satisfied for a infinite range of values.

{\bf The Riemann hypothesis}

Propositions 1 and 2 show that the worst case for the inequality (3), if not satisfied, should occur at a primorial $N_n$. Proposition 3 deals about the distribution of values of $f(N_n)$ at large n. Propositions 1 and 2 are needed in the proof of RH, based on the refined Robin inequality.

Let us first show that RH $\Rightarrow$ (\ref{Dedekind}). 

This is easy because if RH is true, then Robin's inequality (\ref{Robin}) is true, for any $m>5040$, including at primorials $m=N_n$, $m=N_{n+1}$ and so on \cite{Robin84}, despite the result established in proposition 3 that there are infinitely many values of $n$ such that $f(N_{n+1})>f(N_n)$ \footnote {In the first version of this paper, it was expected that for any $p_n$, one would have $\theta(p_n)<p_n$ and as result $f(N_{n+1})<f(N_n)$. But this property is unnecessary for showing that the refined Robin inequality is equivalent to RH.}. Since $N_n$ is free of square, one has $\psi(N_n)=\sigma(N_n)$ so that the refined inequality (\ref{Dedekind}) is satisfied at any $m=N_n$.
This means that if RH is true then, according to proposition 1, (\ref{Dedekind}) is satisfied at $lN_n$, where $N_n<lN_n<N_{n+1}$ and, according to proposition 2, it is also satisfied at any $m$ between consecutive values $lN_n$ and $(l+1)N_n$ of the sequence $A060735$.

The reverse implication (3) $\Rightarrow$ RH is similar to that for the Robin's inequality (theorem 2 in \cite{Robin84}). We observe that there exists an infinity of numbers $n$ such that $\tilde{g}(n)=\frac{\sigma(n)}{n \log \log n}>e^{\gamma}$ and the bound on $\frac{\psi(n)}{n}$ is such that
$$\mbox{for}~n\ge 3,~ \frac{\psi(n)}{n}\le \frac{\sigma(n)}{n}\le e^{\gamma} \log \log n +\frac{0.6482}{\log \log n}.$$

To prove that RH is true, it is sufficient to prove that no exception to the refined Robin's criterion may be found.

At large value of primorials $N_n$, we use Mertens theorem about the density of primes $\log p_n \prod_{k=1}^n (1-\frac{1}{p_k})\sim e^{- \gamma}$, or the equivalent relation \footnote{A better approximation could be obtained  from proposition 9 in \cite{Carella10}, i.e. from $\prod_{p\le x}(1+\frac{1}{p})=\frac{e^{\gamma}}{\zeta (2)} \log x +O(1/ \log x)$.}
\begin{equation}
\frac{1}{\log p_n}\frac{\psi(N_n)}{N_n}\equiv \frac{1}{\log p_n} \prod_{k=1}^n (1+\frac{1}{p_k}) \sim \frac{e^{\gamma}}{\zeta(2)},
\label{Mertens} 
\end{equation}
and the lower bound given p. 206 of \cite{Robin84}

\begin{equation}
\mbox{for}~p_n\ge 20000,~\log \log N_n > \log p_n-\frac{0.123}{\log p_n}.
\label{lowerbound} 
\end{equation}
Using (\ref{Mertens}) and (\ref{lowerbound}), and with $e^{\gamma}(\zeta(2)^{-1}-1)\approx -0.698$, one obtains a lower bound
\begin{equation}
\mbox{for}~p_n\ge 20000,~f(N_n)<-0.698 \log p_n+\frac{0.220}{\log p_n} \sim -6.89.
\label{lbound} 
\end{equation}
For values of $ 2 \le n \le 100000$, computer calculations can be performed. The calculation of $f(N_n)=g(N_n)<0$ is fast using the multiplicative property of $\sigma(n)$, i.e. using $\sigma(N_n)=\prod_{i=1}^n(1+p_i)$. One find a decreasing function $g(N_n)$ that is negative if $n>3$, i.e. $N_n>30$,  as expected (and in agreement with the exceptions in the sets $\mathcal{A}$ and $\mathcal{B}$).

\begin{table}[ht]
\begin{center}
\small
\begin{tabular}{|r|r|r|r|r|r|r|r|}
\hline
$n$ & $3$ & $10$ & $10^2$& $10^3$ & $10^4$& $10^5$\\
\hline
$f(N_n)=g(N_n)$ & $0.22$ & $-1.67$ & $-4.24$& $-6.23$ & -8.06& -9.83 \\
\hline
\end{tabular}
\label{table2}
\caption{The approximate value of the function $f(N_n)=g(N_n)$ versus the number of primes in the primorial $N_n$. The smallest primorial in the table is $N_n=30$ and the highest one is $ N_{100000}\approx 1.9 \times 10^{563920}$.}
\end{center}
\end{table}
According to the calculations illustated in the table II and the bound established in (\ref{lbound}), there are no exceptions to the refined Robin's criterion. Since Robin's criterion has been shown to be equivalent to RH hypothesis, RH may only be true.

\section*{Acknowledgements} The author thanks P. Sol\'e for pointing out his paper \cite{Sole07}, following the presentation of \cite{Planat10} at QuPa meeting in Paris on 09/23/2010. He acknowledges J.~L. Nicolas for helping him to establish the proposition 3. He also thanks Fabio Anselmi for his sustained feedback on this subject.

\end{document}